\documentclass{amsart}
\usepackage{amssymb,latexsym}
\numberwithin{equation}{section}
\newtheorem{theorem}{Theorem}[section]
\newtheorem{proposition}[theorem]{Proposition}

\newtheorem{question}[theorem]{Question}
\newtheorem{corollary}[theorem]{Corollary}
\newtheorem{remark}[theorem]{Remark}
\newtheorem{lemma}[theorem]{Lemma}

\newtheorem{definition}[theorem]{Definition}

\def\proof{\smallskip\noindent {\bf Proof: }}

\newcommand{\ww}{\mathbf{i}}
\begin{document}
\title[A multiplicative property of quantum flag minors II]
{A multiplicative property of quantum flag minors II}
\author{Philippe Caldero}
\address{D\'epartement de Math\'ematiques, Universit\'e Claude Bernard Lyon I,
69622 Villeurbanne Cedex, France}
\email{caldero@igd.univ-lyon1.fr}
\author{Robert Marsh}
\address{Department of Mathematics and Computer Science, University of
Leicester, University Road, Leicester LE1 7RH, England}
\email{rjm25@mcs.le.ac.uk}
\begin{abstract}
Let $U^+$ be the plus part of the quantized enveloping algebra of a simple
Lie algebra and let $\mathcal{B}^*$ be the dual canonical basis of $U^+$.
Let $b$, $b'$ be in $\mathcal{B}^*$ and suppose that one of the two elements
is a $q$-commuting product of quantum flag minors. We show that $b$ and
$b'$ are multiplicative if and only if they $q$-commute.
\end{abstract}
\date{24th January 2003}
\thanks{The authors would like to thank EC grant of the TMR network
``Algebraic Lie Representations", contract no. ERB FMTX-CT97-0100,
for support of the visit of
both authors to Wuppertal in October 2001 and for support of the first
author's visit to Leicester in August 2002, and the University of Leicester
for study leave for the second author in Autumn 2002.}
\maketitle
\tableofcontents
\section{Introduction}
Let $\mathcal{B}^*$ denote the dual of the canonical
basis~\cite{kashiwara1},~\cite[\S14.4.6]{lusztig1} of a quantized enveloping
algebra of a simple Lie algebra $\mathfrak{g}$. 
Two elements of $\mathcal{B}^*$ are said to be multiplicative if their
product also lies in $\mathcal{B}^*$ up to a power of $q$. They are said to
$q$-commute if they commute up to a power of $q$.
The Berenstein-Zelevinsky conjecture states that two elements of
$\mathcal{B}^*$ are multiplicative if and only if they $q$-commute. While
Reineke~\cite[4.5]{reineke1} has shown that if two elements are multiplicative
then they $q$-commute, the conjecture is now known to be false: a
counter-example was provided by Leclerc~\cite{leclerc1}, who showed that for
all but finitely many types of simple Lie algebra, there exist elements
$b\in\mathcal{B}^*$ whose square does not lie in $\mathcal{B}^*$ even up
to a power of $q$ (such elements obviously $q$-commute with themselves).
Such elements are called imaginary; elements of $\mathcal{B}^*$ whose square
does lie in $\mathcal{B}^*$ up to a power of $q$ are known as real.
We consider the following question:
\begin{question}\label{question} Let $b$, $b'$ be in $\mathcal{B}^*$ and
suppose that $b$ is real. Is it the case that $b$ and $b'$ are multiplicative
if and only if they $q$-commute?
\end{question}
Leclerc has made a conjecture~\cite[Conjecture 1]{leclerc1} concerning the
expansion of the product of two dual canonical basis elements, and has
remarked that if true, it would imply that the answer to this question is yes
in general.
In this paper, we prove that the answer to this question is yes if one of the
two elements involved is a $q$-commuting product of quantum flag minors
(in type $A$). We note that such products are known to be real~\cite{lnt1}.
A key part of our proof involves results concerning the piecewise-linear
repa\-ra\-me\-tri\-za\-tion function $R_{\ww}^{\ww'}$ associated by Lusztig
to a pair $\ww$, $\ww'$ of reduced decompositions for the longest word $w_0$ in
the Weyl group of $\mathfrak{g}$. Lusztig has defined the canonical basis
via the bases of Poincar\'e-Birkhoff-Witt type associated to such
reduced decompositions, and such reparametrization functions arise from taking
two of the Lusztig parametrizations of the canonical basis. We show that
these functions share some of the properties known to be possessed by the
reparametrization functions arising from string parametrizations for the
canonical basis~\cite{berensteinzelevinsky1}. In particular we show that
if two dual canonical basis elements $q$-commute, then their PBW
parametrizations (with respect to a fixed reduced decomposition $\ww$) lie
in a single PBW $\ww$-linearity domain; in other words, if their
parameters are $\mathbf{m}$ and $\mathbf{m}'$, then
$R_{\ww}^{\ww'}(\mathbf{m}+\mathbf{m}')=R_{\ww}^{\ww'}(\mathbf{m})+
R_{\ww}^{\ww'}(\mathbf{m}')$ for all reduced decompositions $\ww'$ for $w_0$.
The corresponding result for string parametrizations is already
known~\cite[2.9]{berensteinzelevinsky1}. Such behaviour hints at an
explanation for the compatibility of examples in which the canonical
basis has been computed
explicitly~\cite{cartermarsh1},~\cite{lusztig3},~\cite{xi1}
with respect to linearity domains.
As a consequence, we obtain that the set of PBW $\ww$-linearity domains forms
a  fan in $\mathbb{R}^N$, where $N$ is the length of $w_0$.\par
This enables us to prove that the answer to Question~\ref{question} is
yes in type A$_n$ when $b$ is a $q$-commuting product of quantum flag minors.
This generalizes a theorem of ~\cite{caldero2}.
\par
In the paper, the Lie algebra $\mathfrak{g}$ is supposed to be of type A$_n$. 
Note that results in section 1-4 can be easily generalized to all simply-laced 
types. Note also that Lemma~\ref{flagminors} is only true in type A$_n$. Hence, 
Theorem~\ref{main} needs this assumption. 
\section{Background and notation}
We use the set-up of~\cite{caldero1}.
Let $\mathfrak{g}=sl_{n+1}(\mathbb{C})$ denote the simple Lie algebra of type
$A_n$. Let $\mathfrak{h}$ be a Cartan subalgebra and let
$\mathfrak{g}=\mathfrak{n}^-\oplus \mathfrak{h}\oplus \mathfrak{n}^+$
be a compatible triangular decomposition. Let
$\alpha_1,\alpha_2,\ldots ,\alpha_n$ be the corresponding simple roots,
and $\Delta^+$ the corresponding set of positive roots.
Let $W$ be the Weyl group associated to the root system, $P$ be the weight 
lattice generated by the fundamental weights $\varpi_i$, $1\leq i\leq n$, and 
$\langle\,,\,\rangle$ be the $W$-invariant form on $P$.
Let $U=U_q(\mathfrak{g})$ be the simply connected Drinfel'd-Jimbo quantized
enveloping algebra over $\mathbb{Q}(q)$ associated to $\mathfrak{g}$
as defined in~\cite{joseph1}. Let $U^-$, $U^0$ and $U^+$ be the subalgebras
associated to the sub-Lie algebras
$\mathfrak{n}^-$, $\mathfrak{h}$ and $\mathfrak{n}^+$
respectively; we have the triangular decomposition
$U\cong U^-\otimes U^0\otimes U^+$.
The subalgebra $U^+$ is generated over $\mathbb{Q}(q)$ by canonical
generators $E_1,E_2,\ldots E_n$, subject to the quantized Serre relations;
the subalgebra $U^-$ is isomorphic to $U^+$, with corresponding generators
$F_1,F_2,\ldots F_n$, and the subalgebra $U^0$ is isomorphic to
$\mathbb{Q}(q)[P]$, the element corresponding to $\lambda\in P$ being
denoted by $K_{\lambda}$. Let $\mathfrak{b}^+=\mathfrak{h}\oplus \mathfrak{n}^+$ and
$\mathfrak{b}^-=\mathfrak{h}\oplus \mathfrak{n}^-$ and we define
$U(\mathfrak{b}^+)=U^0U^+$ and $U(\mathfrak{b}^-)=U^-U^0$.
The algebra $U$ is a Hopf algebra with comultiplication $\Delta$,
antipode $S$ and augmentation $\varepsilon$ given by:
$$\Delta(E_i)=E_i\otimes 1+K_{\alpha_i}\otimes E_i,\ \ 
\Delta(F_i)=F_i\otimes K_{-\alpha_i}+1\otimes F_i,\ \ 
\Delta(K_{\lambda})=K_{\lambda}\otimes K_{\lambda},$$
$$S(E_i)=-K_{-2\alpha_i}E_i,\ \ 
S(F_i)=-F_iK_{2\alpha_i},\ \ 
S(K_{\lambda})=K_{-\lambda},$$
$$\varepsilon(E_i)=\varepsilon(F_i)=0,\ \ 
\varepsilon(K_{\lambda})=1.$$
The root lattice $Q$ is defined to be $Q=\sum_i \mathbb{Z}\alpha_i$,
with $Q^+=\sum_i \mathbb{Z}_{\geq 0}\alpha_i$.
Recall that if $\alpha=\sum_i m_i\alpha_i\in Q^+$, 
then an element in the subspace of $U$ generated
by $\{E_{i_1}^{n_1}E_{i_2}^{n_2}\cdots E_{i_k}^{n_k}\,:\,
n_1\alpha_{i_1}+n_2\alpha_{i_2}+\cdots +n_k\alpha_{i_k}=\alpha\}$
(respectively,
$\{F_{i_1}^{n_1}F_{i_2}^{n_2}\cdots F_{i_k}^{n_k}\,:\,
n_1\alpha_{i_1}+n_2\alpha_{i_2}+\cdots +n_k\alpha_{i_k}=\alpha\}$),
is said to have weight $\alpha$ (respectively, $-\alpha$).
For all $\alpha\in Q^+$, let $U^+_{\alpha}$ (respectively, $U^-_{-\alpha}$)
denote the subspace of $U_+$ (respectively, $U^-$) consisting of elements of
weight $\alpha$ (respectively, $-\alpha$). We write the weight of an element
$X$ as $wt(X)$ (if it exists), and $tr(X)$ (or $tr(wt(X))$) for the sum
$\sum_i m_i$ if $wt(X)=\sum_i m_i\alpha_i$.
For $u\in U$, set $\Delta(u)=u_{(1)}\otimes u_{(2)}\in U\otimes U$.
There exists a unique bilinear form $(\,,\,)$~\cite{rosso1},~\cite{tanisaki1}
on $U(\mathfrak{b}^+)\times U(\mathfrak{b}^-)$ satisfying
$$(E_i,F_i)=\delta_{ij}(1-q^2)^{-1},$$
$$(u^+,u_1^-u_2^-)=(\Delta(u^+),u_1^-\otimes u_2^-),\ \ \ \ \ \ \ \ 
u^+\in U(\mathfrak{b}^+), u_1^-,u_2^-\in U(\mathfrak{b}^-),$$
$$(u_1^+u_2^+,u^-)=(u_2^+\otimes u_1^+,\Delta(u^-)),\ \ \ \ \ \ \ \ 
u^-\in U(\mathfrak{b}^-), u_1^+,u_2^+\in U(\mathfrak{b}^+),$$
$$(K_{\lambda},K_{\mu})=q^{-(\lambda,\mu)},\ \ (K_{\lambda},F_i)=0,\ \ 
(E_i,K_{\lambda})=0,\ \ \ \ \ \ \ \ \lambda,\mu\in P.$$
The form $(\,,\,)$ is nondegenerate on $U^+_{\alpha}\otimes U^-_{-\alpha}$
for all $\alpha\in Q_+$. Since $U^+$ and $U^-$ are isomorphic algebras
(with isomorphism preserving their weight spaces), we can naturally identify
$U^-_{-\alpha}$ with $U^+_{\alpha}$, so for each element $u\in U^+_{\alpha}$
there is a corresponding element $u^*\in U^+_{\alpha}$ (corresponding to
$u$ using the form). 
Let $w_0$ be the longest element of $W$;
denote by $R(w_0)$ the set of all reduced decompositions for $w_0$.
Fix a reduced decomposition $\ww=(i_1,i_2,\ldots ,i_N)$ of $w_0$, and for
$1\leq t\leq N$, let
$\beta_t=s_{i_1}s_{i_2}\cdots s_{i_{t-1}}(\alpha_{i_t})$; we get an ordering
$\beta_1<\beta_2<\cdots <\beta_N$ of $\Delta^+$. For $1\leq i\leq n$, let
$T_i$ denote the Lusztig braid
automorphism of $U^+$~\cite[37.1.3]{lusztig1},~\cite{saito1}
associated to $i$. For
$1\leq t\leq N$ let $E_{\beta_t}=T_{i_1}T_{i_2}\cdots T_{i_{t-1}}(E_{i_t})$.
The Poincar\'e-Birkhoff-Witt basis of $U^+$ is the basis
$$B_{\ww}=\{E(\mathbf{m})\,:\,\mathbf{m}\in \mathbb{Z}_{\geq 0}^N\},$$
where
$$E(\mathbf{m})=E_{\ww}(\mathbf{m})=\prod_{t=1}^N
\frac{1}{[m_t]!}E_{\beta_t}^{m_t},$$
with the product taken in the ordering given above. Here
$[m]!=[m][m-1]\cdots [1]$, where $[m]=\frac{q^m-q^{-m}}{q-q^{-1}}$.
By~\cite{levendorskiisoibelman1} we have that
$$E(\mathbf{m})^*=\prod_{t=1}^N\psi_{m_t}(q^2)E(\mathbf{m}),$$
where $\psi_m(z)=\prod_{k=1}^m (1-z^k)$.
Let $\mathcal{L}$ be the sub-$\mathbb{Z}[q]$-lattice of $U^+$ generated
by $B_{\ww}$, and let $\mathcal{L}^*$ be the sub-$\mathbb{Z}[q]$-lattice of
$U^+$ generated by $B_{\ww}^*=\{E_{\ww}(\mathbf{m})^*\,:\,\mathbf{m}\in
\mathbb{Z}_{\geq 0}^N\}$. Lusztig has shown that both $\mathcal{L}$ and
$\mathcal{L}^*$ are independent of the choice of reduced decomposition $\ww$.
Let $\eta$ be the $\mathbb{Q}$-algebra automorphism of $U$ fixing the
generators $E_i$ and $F_i$, with $\eta(K_{\lambda})=K_{-\lambda}$ and
$\eta(q)=q^{-1}$. Lusztig~\cite{lusztig2} also shows:
\begin{theorem} (Lusztig)
Fix a reduced decomposition $\ww$ of $w_0$. Then, for each $\mathbf{m}\in
\mathbb{Z}_{\geq 0}^N$, there is a unique element $B(\mathbf{m})=
B_{\ww}(\mathbf{m})\in U^+$ such that $\eta(B(\mathbf{m}))=B(\mathbf{m})$
and $B(\mathbf{m})\in E(\mathbf{m})+q\mathcal{L}$. The set
$\mathcal{B}=\{B(\mathbf{m})\,:\,\mathbf{m}\in \mathbb{Z}_{\geq 0}\}$
is a basis (called the canonical basis) of $U^+$ which does not depend
on $\ww$.
\end{theorem}
We'll call the parametrization $\mathbb{Z}_{\geq 0}^N\rightarrow \mathcal{B}$,
$\mathbf{m}\mapsto B(\mathbf{m})$,
Lusztig's parametrization of $\mathcal{B}$ (arising from $\ww$).
If $b\in B$, we denote by $L_{\ww}(b)$ its Lusztig parameter (so $L_{\ww}$
is the inverse of $B_{\ww}$).
For any $\ww,\ww'\in R(w_0)$, we also have Lusztig's piecewise-linear
reparametrization function $R_{\ww}^{\ww'}=L_{\ww'}L_{\ww}^{-1}$ (which can
be regarded as a function from $\mathbb{R}^{N}\rightarrow \mathbb{R}^N$
via the same formula in terms of coordinates).
The canonical basis was discovered independently by
Kashiwara~\cite{kashiwara1}, who called
it the global crystal basis. Let $\mathcal{B}^*=\{b^*\,:\,b\in\mathcal{B}\}$
denote the dual canonical basis of $U^+$. For $b\in \mathcal{B}$ we define
$L_{\ww}(b^*)$ to be $L_{\ww}(b)$. Let $\sigma$ be the antihomomorphism
of $U$ (as $\mathbb{Q}(q)$-algebra) taking $E_i$ to $E_i$, $F_i$ to $F_i$
and $K_{\lambda}$ to $K_{-\lambda}$. As in~\cite[Proposition 16]{lnt1} we have:
\begin{proposition} \label{dualclass}
Fix a reduced decomposition $\ww$ of $w_0$. Then, for each $\mathbf{m}\in
\mathbb{Z}_{\geq 0}^N$, the element $B(\mathbf{m})^*$ is the unique element
$X$ of $U^+$ with weight $\sum_i m_i\alpha_i$ such that
\begin{equation} \label{dualclasseq}
\eta(X)=(-1)^{tr(X)}q^{-\langle wt(X),wt(X)\rangle/2}q_X^{-1}\sigma(X),\ 
X\in E(\mathbf{m})^*+q\mathcal{L}^*,
\end{equation}
where $q_X=\prod_{i}q^{m_i}$, $wt(X)=\sum_i m_i\alpha_i$.
\end{proposition}
Fix a reduced decomposition $\ww$ in $R(w_0)$. Let $\mathbf{e}_k$,
$1\leq k\leq N$, be the canonical generators of $\mathbb{Z}_{\geq 0}^N$.
We can define an ordering on the semigroup $\mathbb{Z}_{\geq 0}^N$ associated
to $\ww$ in the following way. The ordering $\prec_{\ww}$ is generated by
$\mathbf{m}\leq \mathbf{e}_k+\mathbf{e}_{k'},\;
1\leq k<k'\leq N\Leftrightarrow E_{\ww}(\mathbf{m})$ is a term of the
PBW decomposition of
$E_{\beta_{k}}E_{\beta_{k'}}-q^{(\beta_{k'},\beta_{k})}E_{\beta_{k'}}
E_{\beta_{k}}$.
It will be denoted by $\prec$ if no confusion occurs.
\begin{proposition} \label{ordering}
Fix $\ww$ in $R(w_0)$. Then, for all $\bf{m}$, $\bf{n}$ in
$\mathbb{Z}_{\geq 0}^N$:\par\noindent
\item{(i)} if $\mathbf{m}\prec_{\ww}\mathbf{n}$ then $\mathbf{m}$ is lower
than $\mathbf{n}$ for the lexicographical ordering,\par\noindent
\item{(ii)} $B_{\ww}(\mathbf{m})=E_{\ww}(\mathbf{m})+\sum_{\mathbf{m}
\prec\mathbf{n}}d_{\mathbf{m}}^{\mathbf{n}}E_{\ww}(\mathbf{n})$,
$d_{\mathbf{m}}^{\mathbf{n}}\in q\mathbb{Z}[q]$,\par\noindent
\item{(iii)} $B_{\ww}(\mathbf{m})^*=E_{\ww}(\mathbf{m})^*+q\sum_{\mathbf{n}
\prec\mathbf{m}}c_{\mathbf{m}}^{\mathbf{n}}E_{\ww}(\mathbf{n})^*$,
$c_{\mathbf{m}}^{\mathbf{n}}\in \mathbb{Z}[q]$. 
\end{proposition}
\proof
See~\cite[2.1]{caldero2}.~$\Box$
\par
Let $q^{\mathbb{Z}}\mathcal{B}^*$ denote the set $\{q^n b^*\,:\,b\in 
\mathcal{B},n\in\mathbb{Z} \}$. We therefore have:
\begin{corollary} \label{dualclassq}
Let $\alpha\in Q^+$, and suppose $u\in U^+_{\alpha}$. Fix $\ww\in R(w_0)$.
Then $u\in q^{\mathbb{Z}}\mathcal{B}^*$ if and only if \\
(i) There is $f\in \mathbb{Z}[q,q^{-1}]$ such that
$\sigma\eta(u)=fu$, and \\
(ii) $u=q^k(E(\mathbf{m})^*+q\sum_{\mathbf{n}<\mathbf{m}}
c_{\mathbf{m}}^{\mathbf{n}} E_{\ww}(\mathbf{n})^*)$,
where $k\in \mathbb{Z}$, $c_{\mathbf{m}}^{\mathbf{n}}\in \mathbb{Z}[q]$
and $<$ denotes the lexicographic ordering.
\end{corollary}
\proof It is immediate from Proposition~\ref{dualclass} that if $u\in
q^{\mathbb{Z}}\mathcal{B}^*$ then it satisfies (i) and (ii).
If $k=0$ in part (ii), this ensures that, in the expansion of $u$ in
terms of the dual canonical basis, $B_{\ww}(\mathbf{m})^*$ occurs with
coefficient $1$ (see Proposition~\ref{ordering}). It follows that the
eigenvalue $f$ in part (i) must be the same as that appearing
in~(\ref{dualclasseq}), and it follows from Proposition~\ref{dualclass}
that $u\in \mathcal{B}^*$. The result for arbitrary $k$ follows.~$\Box$
\begin{remark} \rm
We note that, from the proof, we can see that $u\in \mathcal{B}^*$ if and
only if $u$ satisfies (i) and (ii) with $k=0$.
\end{remark}
Two elements of the dual canonical basis $\mathcal{B}^*$ are
said to be multiplicative if their product also lies in the dual canonical
basis up to a power of $q$. They are said to $q$-commute if they commute
up to a power of $q$. The following Corollary of
Proposition~\ref{dualclass} is due to Reineke~\cite[4.5]{reineke1}. 
\begin{corollary}
If two elements of the dual canonical basis are multiplicative then they
$q$-commute.
\end{corollary}
\begin{definition} \label{action}
There is a natural action of $U^+$ on itself~\cite{berensteinzelevinsky1},
in which each generator $E_i$ of $U^+$ acts on $U^+$ as a $q$-differential
operator $\delta_i$. This is the unique action of $U^+$ on itself satisfying
the following properties: \\
(a) (Homogeneity) If $E\in U^+_{\alpha}$, $x\in U^+_{\gamma}$, then
$E(x)\in U^+_{\gamma-\alpha}$. \\
(b) (Leibnitz formula)
$$\delta_i(xy)=\delta_i(x)y+q^{-(\gamma,\alpha_i)}x\delta_i(y),\mbox{\ for\ 
}x\in
U^+_{\gamma},\ y\in U^+.$$
(c) (Normalization) $\delta_i(E_j)=(1-q^2)^{-1}\delta_{ij}$ for $i,j=1,2,\ldots 
,n$.
\end{definition}
We remark that formula (b) implies that for all $i$ and for
$r\in\mathbb{Z}_{\geq 0}$, we have 
$\delta_i(E_i^{(r)})=q^{-r+1}(1-q^2)^{-1}E_i^{(r-1)}$,
which is easily checked by induction on $r$ ($E_i^{(0)}$ is interpreted as
$1$ and $E_i^{(-1)}$ as zero).
For $u\in U^+$, set $\varphi_i(u)=\max\{r,\delta_i^r(u)\not=0\}$,
and for $r\in\mathbb{Z}_{\geq 0}$, let $\delta_i^{(r)}$ denote the divided
power $\frac{\delta_i^r}{[r]!}$.
Define $\delta_i^{(max)}(u):=
\delta_i^{(\varphi_i(b^*))}(u)$. It is known that
for $b\in\mathcal{B}$, $\delta_i^{(max)}(b^*)\in \mathcal{B}^*$
(see~\cite[\S 1]{berensteinzelevinsky1}).
If $\ww\in R(w_0)$, set $a_1=\varphi_{i_1}(b^*)$,
$a_2=\varphi_{i_2}(\delta_{i_1}^{(a_1)}(b^*)),\ldots ,a_N=
\varphi_{i_N}(\delta_{i_{N-1}}^{(a_{N-1})}\cdots \delta_{i_1}^{(a_1)}(b^*))$;
then $(a_1,a_2,\ldots ,a_N)$ is known as the string of $b$ (or $b^*$) in
direction $\ww$~\cite{berensteinzelevinsky1}; this coincides with the string
of $b$ arising from Kashiwara's approach to $\mathcal{B}$
(see~\cite{kashiwara2} and the end of Section~$2$
in~\cite{nakashimazelevinsky1}).
%***************************************************************************
%***************************************************************************
\section{PBW-strings}
In this section we will show how the Lusztig parametrization of a dual canonical
basis element can also be regarded as a string (in a similar sense to the
above). We define new operators $\Delta_i$, $i=1,2,\ldots ,n$, depending
on a choice of reduced decomposition for $w_0$, which play
the role of the operators $\delta_i^{(max)}$ for the PBW parametrization.
We first of all note that $\delta_i^{(max)}$ is a well-defined operator on
all of $U^+$.
\begin{definition} \label{Delta}
Let $\widetilde{w}=s_{i_1}s_{i_2}\cdots s_{i_t}$ be a reduced decomposition for
$w\in W$. Let $\Delta_1=\delta_{i_1}^{(max)}$,
$\Delta_2=T_{i_1}\delta_{i_2}^{(max)}T_{i_1}^{-1},\ldots ,
\Delta_t=(T_{i_1}\cdots T_{i_{t-1}})\delta_{i_t}^{(max)}
(T_{i_1}\cdots T_{i_{t-1}})^{-1}$, operators on $U^+$ (with codomain $U$).
\end{definition}
We will use the following:
\begin{lemma} \label{braidimage}
(Saito) Let $i\in \{1,2,\ldots ,n\}$. Then $T_i(U^+)\cap U^+=\ker \delta_i$.
\end{lemma}
\proof See~\cite{saito1}.~$\Box$
\begin{lemma} \label{deltaPBW}
Let $\ww\in R(w_0)$, and $\mathbf{m}=(m_1,m_2,\ldots ,m_N)\in
\mathbb{Z}_{\geq 0}^N$. Then
$$\delta_{i_1}^{(max)}(E_{\ww}(\mathbf{m})^*)=
E_{\ww}(0,m_2,m_3,\ldots ,m_N)^*.$$
\end{lemma}
\proof Let $\gamma_{\mathbf{m}}=
\prod_{t=1}^N\psi_{m_t}(q^2),$
so that $E_{\ww}(\mathbf{m})^*=
\gamma_{\mathbf{m}}E_{\ww}(\mathbf{m})$.
Then $E_{\ww}^{\mathbf{m}}=E_{i_1}^{(m_1)}T_{i_1}(y)$, where
$y=E_{i_2}^{(m_2)}T_{i_2}(E_{i_3}^{(m_3)})\cdots T_{i_2}T_{i_3}\cdots
T_{i_{N-1}}(E_{i_N}^{(m_N)})\in U^+$. Using Lemma \ref{braidimage}, the
formula for $\gamma_{\mathbf{m}}$ and the fact
that $\delta_{i_1}(E_{i_1}^{(m)})=q^{-m+1}(1-q^2)^{-1}E_{i_1}^{(m-1)}$ the 
result
follows.~$\Box$
\par
We can now prove a Lemma giving basic properties of the Lusztig parametrization
with respect to the operators $\Delta_i$.
\begin{lemma} \label{PBWstring}
Let $b\in\mathcal{B}$, and suppose $\ww\in R(w_0)$
and $\Delta_1,\Delta_2,\ldots ,\Delta_N$ are as above. Let
$L_{\ww}(b)=(m_1,m_2,\ldots ,m_N)\in \mathbb{Z}_{\geq 0}^N$.
Then, for $k=1,2,\ldots ,N$, we have: \\
(i) $\Delta_k\Delta_{k-1}\cdots \Delta_1(b^*)\in \mathcal{B}^*$,
\\
(ii)
$\varphi_{i_k}((T_{i_1}\cdots T_{i_{k-1}})^{-1}\Delta_{k-1}\cdots 
\Delta_1(b^*))=m_k$, \\
(iii) $L_{\ww}(\Delta_k\Delta_{k-1}\cdots \Delta_1(b^*))=(0,0,\ldots, 0,m_{k+1},
m_{k+2},\ldots ,m_N)$.
\end{lemma}
\proof
Since $\delta_{i_1}^{(max)}$ preserves $\mathcal{B}^*$,
we see that $\Delta_1(b^*)$ lies in $\mathcal{B}^*$,
so (i) holds for $k=1$, but we need more precise information.
Since $b^*\in \mathcal{B}^*$, we know by~\cite{caldero2}
(see also Proposition~\ref{ordering}) that
$$b^*=E_{\ww}(\mathbf{m})^*+qx_0,$$
where $x_0$ is a $\mathbb{Z}[q]$-linear combination of dual PBW-basis elements
$E_{\ww}(\mathbf{n})^*$ with $\mathbf{n}<\mathbf{m}$ (and
in particular, satisfying $n_1\leq m_1$). Here $<$ denotes the lexicographic
ordering. By Lemma~\ref{deltaPBW}, we see that
$$\delta_{i_1}^{(max)}(b^*)=E_{\ww}(0,m_2,m_3,\ldots ,m_N)^*+qx_1,$$
where $x_1$ is a $\mathbb{Z}[q]$-linear combination of elements of the form
$E_{\ww}(\mathbf{n})^*$ with
$\mathbf{n}=(0,n_2,n_3,\ldots ,n_N)<\mathbf{m}$ for all $\mathbf{n}$ occurring
in the sum.
By Corollary~\ref{dualclassq}, we have that
$L_{\ww}(\Delta_1(b^*))=(0,m_2,m_3,\ldots ,m_N)$, so (ii) holds for $k=1$.
It is also clear that $\varphi_{i_1}(b^*)=m_1$, so (iii) holds for $k=1$.
We thus see that (i), (ii) and (iii) all hold for $k=1$.
\par
>From the above form for $\Delta_1(b^*)$, we have that
$$T_{i_1}^{-1}(\Delta_1(b^*))=E_{\ww'}(m_2,m_3,\ldots ,m_N,0)^*+qy_1,$$
where $\ww'=s_{i_2}s_{i_3}\cdots s_{i_N}s_{i_1^*}$,
$i_1^*$ is the index of the simple root $-w(\alpha_{i_1})$ (the
Chevalley automorphism applied to $i_1$), and $y_1$ is a
$\mathbb{Z}[q]$-linear combination
of dual PBW-basis elements $E_{\ww'}(\mathbf{n})^*$ with
$\mathbf{n}=(n_2,n_3,\ldots ,n_N,0)<(m_2,m_3,\ldots ,m_N,0)$ and therefore
satisfies Corollary~\ref{dualclassq}(ii) (with $k=0$).
Since for all $i$, $T_i^{-1}\sigma\eta$ and $\sigma\eta T_i^{-1}$
differ on weight spaces by plus or minus a power of $q$,
it follows that $T_{i_1}^{-1}(\Delta_1(b^*))$ satisfies
Corollary~\ref{dualclassq}(i), and therefore (see the remark after 
Corollary~\ref{dualclassq}(i)) lies in
$\mathcal{B}^*$.
Since $\delta_{i_2}^{(max)}$ preserves $\mathcal{B}^*$, we thus obtain that
$\delta_{i_2}^{(max)}T_{i_1}^{-1}\Delta_1(b^*)\in \mathcal{B}^*$.
Arguing as above, with $i_2$ playing the role of $i_1$, we obtain that
$$\delta_{i_2}^{(max)}T_{i_1}^{-1}\Delta_1(b^*)=
E_{\ww'}(0,m_3,m_4,\ldots ,m_N,0)^*+qy_2,$$
where $y_2$ is a $\mathbb{Z}[q]$-linear combination of dual PBW-basis elements
of the form $E_{\ww'}(\mathbf{n})^*$ with
$\mathbf{n}=(0,n_3,n_4,\ldots ,n_N,0)<(0,m_3,m_4,\ldots ,m_N,0)$.
We also obtain that
$\varphi_{i_2}(T_{i_1}^{-1}\Delta_1(b^*))=m_2$ (part (ii) for $k=2$).
Applying $T_{i_1}$ we obtain that
$$\Delta_2\Delta_1(b^*)=T_{i_1}\delta_{i_2}^{(max)}T_{i_1}^{-1}\Delta_1(b^*)=
E_{\ww}(0,0,m_3,m_4,\ldots ,m_N)^*+qx_2,$$ where $x_2$ is a
$\mathbb{Z}[q]$-linear combination of dual PBW-basis elements of the form
$E_{\ww}(\mathbf{n})^*$ with $\mathbf{n}=(0,0,n_3,n_4,\ldots ,n_N)<
(0,0,m_3,m_4,\ldots ,m_N)$. Arguing as above for $T_{i_1}^{-1}$ and
applying Corollary~\ref{dualclassq}, we obtain
that $\Delta_2\Delta_1(b^*)\in \mathcal{B}^*$ and that
$L_{\ww}(\Delta_2\Delta_1(b^*))=(0,0,m_3,m_4,\ldots ,m_N)$.
We thus see that (i) and (iii) hold for $k=2$. It is now clear that an
inductive argument gives (i),(ii) and (iii) for $k=1,2,\ldots ,N$.~$\Box$
\par
Given a reduced decomposition $\widetilde{w}=s_{i_1}s_{i_2}\cdots s_{i_t}$
for an element $w\in W$,
we define the PBW-string $L_{\widetilde{w}}(b^*)$ of a dual canonical
basis element $b^*$ in direction $\widetilde{w}$ as follows. Let
$\ww=s_{i_1}s_{i_2}\cdots s_{i_t}s_{i_{t+1}}\cdots s_{i_N}$ be
any completion of $\widetilde{w}$ to a reduced decomposition for $w_0$.
Then let $L_{\widetilde{w}}(b^*)=(m_1,m_2,\ldots ,m_t)$ where
$L_{\ww}(b^*)=(m_1,m_2,\ldots ,m_N)$. It is clear from Lemma~\ref{PBWstring}
that this is well-defined.
\par
We make definitions for PBW-strings in the same way as Berenstein and
Zelevinsky~\cite{berensteinzelevinsky1} for usual strings.
\begin{definition}
Let $\ww\in R(w_0)$. A PBW $\ww$-wall
is defined to be a hyperplane in $\mathbb{R}^N$ given by the equation
$a_k=a_{k+2}$ for some index $k$ such that $i_k=i_{k+2}=i_{k+1}\pm 1$.
Let $\mathbf{m}\in \mathbb{R}_{>0}^N$. We say that $\mathbf{m}$ is
PBW $\ww$-regular if for every $\ww'\in R(w_0)$
the point $R_{\ww}^{\ww'}(\mathbf{m})$ does not lie on any PBW $\ww'$-wall.
We define the PBW $\ww$-linearity domains to be the closures of the connected
components of the set of PBW $\ww$-regular points.
\end{definition}
As in the string case~\cite[2.8]{berensteinzelevinsky1},
we have the following justification for the terminology; the proof is
basically the same, as $R_{\ww}^{\ww'}$ is very similar to the
string reparametrization function.
\begin{proposition} \label{linearitydomain}
Every PBW $\ww$-linearity domain is a polyhedral convex cone in
$\mathbb{R}_{\geq 0}^N$. Two points $\mathbf{m},\mathbf{m}'$ in
$\mathbb{R}_{\geq 0}^N$ lie in a single PBW $\ww$-linearity domain
if and only if
$$R_{\ww}^{\ww'}(\mathbf{m}+\mathbf{m}')=R_{\ww}^{\ww'}(\mathbf{m})+
R_{\ww}^{\ww'}(\mathbf{m}')$$
for every $\ww'\in R(w_0)$.
\end{proposition}
We can now prove the analogue of~\cite[2.9]{berensteinzelevinsky1}:
\begin{theorem}\label{analogue}
Let $b^*,b'^*$ be elements of the dual canonical basis that $q$-commute. Then,
for every $\ww\in R(w_0)$, the PBW strings $\mathbf{m}=L_{\ww}(b^*)$ and
$\mathbf{m}'=L_{\ww}(b'^*)$ belong to a single $\ww$-linearity domain.
\end{theorem}
\proof
We follow the proof of Berenstein and Zelevinsky~\cite{berensteinzelevinsky1};
however there will be some differences, so we include the details.
We know by Proposition~\ref{linearitydomain} that it is enough to show that,
for every
$\ww'\in R(w_0)$, 
$R_{\ww}^{\ww'}(\mathbf{m}+\mathbf{m}')=R_{\ww}^{\ww'}(\mathbf{m})+
R_{\ww}^{\ww'}(\mathbf{m}')$.
Thus, it is enough to show that $\mathbf{m}$ and $\mathbf{m}'$ are not separated 
by any
PBW $\ww$-wall.
Suppose that a $\ww$-wall corresponds to the move
$$\ww=s_{i_1}\cdots s_{i_t}s_is_js_is_{i_{t+4}}\cdots s_{i_N}
\mapsto s_{i_1}\cdots s_{i_t}s_js_is_js_{i_{t+4}}\cdots s_{i_N}=\ww'.$$
By~\cite[3.6]{berensteinzelevinsky1} and the fact that the $T_i$ are algebra
automorphisms, the elements
$(T_{i_k}T_{i_{k-1}}\cdots T_{i_1})^{-1}$
$\Delta_t\Delta_{t-1}\cdots\Delta_1(b^*)$ and
$(T_{i_k}T_{i_{k-1}}\cdots T_{i_1})^{-1}$
$\Delta_t\Delta_{t-1}\cdots\Delta_1(b'^*)$ $q$-commute.
Replacing $b^*$ and $b'^*$ by
$(T_{i_k}T_{i_{k-1}}\cdots T_{i_1})^{-1}$
$\Delta_t\Delta_{t-1}\cdots\Delta_1(b^*)$ and
$(T_{i_k}T_{i_{k-1}}\cdots T_{i_1})^{-1}$
$\Delta_t\Delta_{t-1}\cdots\Delta_1(b'^*)$,
we can assume that $t=0$ (using Lemma \ref{PBWstring}).
Let $L_{iji}=L_{s_is_js_i}(b^*)=(m_1,m_2,m_3)$ and
$L'_{iji}=L_{s_is_js_i}(b'^*)=(m'_1,m'_2,m'_3)$
be the PBW-strings of $b^*$ and $b'^*$ in direction $s_is_js_i$.
It is sufficient to show that the integers $m_1-m_3$ and $m'_1-m'_3$ are
of the same sign.
Now let $\Delta_1,\Delta_2,\Delta_3$ be the operators on $\mathcal{B}^*$
associated to the reduced decomposition $s_is_js_i$ as in 
Definition~\ref{Delta}.
Let $b_0^*=\Delta_3\Delta_2\Delta_1(b^*)$, and let
$b'^*_0=\Delta_3\Delta_2\Delta_1(b'^*)$. Suppose that $b'^*b^*=q^nb^*b'^*$
and
$b'^*_0b_0^*=q^{n_0}b_0^*b'^*_0$ for integers $n,n_0$.
Then we have
\begin{eqnarray*}
\Delta_3\Delta_2\Delta_1(b^*b'^*) & = &
T_iT_j\delta_i^{(max)}T_j^{-1}T_i^{-1}\delta_j^{(max)}T_i^{-1}\delta_i^{(max)}(b
b')
\\
& = & q^{r_1}\Delta_3\Delta_2\Delta_1(b^*)\Delta_3\Delta_2\Delta_1(b'^*) \\
& = & q^{r_1}b_0^*b'^*_0,
\end{eqnarray*}
where
$$r_1=\Phi_{i,\gamma}(m_1,m_1')+\Phi_{j,s_i(\gamma-m_1\alpha_i)}(m_2,m'_2)+
\Phi_{i,s_js_i(\gamma-m_1\alpha_i)-m_2\alpha_j},$$
and $\gamma$ is the degree of $b^*$. Here
$\Phi_{k,\mu}(n,m)=nm-(\mu,m\alpha_k)=m(n-(\mu,\alpha_k))$ is defined as
in~\cite[Proposition 3.1]{berensteinzelevinsky1} and we are
using~\cite[(3.6)]{berensteinzelevinsky1}. We also use the fact that
for all $\alpha\in Q^+$,
$T_i^{\pm 1}(U^+_{\alpha})\cap U^+\subset U^+_{s_i(\alpha)}$.
Expanding and simplifying, we obtain
\begin{eqnarray*}
r_1 & = & m_1m'_1+m_2m'_2+m_3m'_3+m_1m'_2+
m_2m'_3-m_1m'_3- \\
& & m'_1(\gamma,\alpha_i)-m'_2(\gamma,\alpha_i+\alpha_j)-m_3(\gamma,\alpha_j).
\end{eqnarray*}
Similarly, we obtain that $\Delta_3\Delta_2\Delta_1(b^*b'^*)=
q^{r_2}\Delta_3\Delta_2\Delta_1(b^*)\Delta_3\Delta_2\Delta_1(b'^*)$,
where
\begin{eqnarray*}
r_2 & = & m'_1m_1+m'_2m_2+m'_3m_3+m'_1m_2+m'_2m_3-m'_1m_3- \\
& & m_1(\gamma',\alpha_i)-m_2(\gamma',\alpha_i+\alpha_j)-
m_3(\gamma',\alpha_j),
\end{eqnarray*}
and $\gamma'$ is the weight of $b'^*$. We thus have:
\begin{equation}\label{difference1}
n_0-n=r_1-r_2.
\end{equation}
Let $L_{jij}=L_{s_js_is_j}(b^*)=(n_1,n_2,n_3)$ and
$L'_{jij}=L_{s_js_is_j}(b'^*)=(n'_1,n'_2,n'_3)$
be the PBW-strings of $b^*$ and $b'^*$ in direction $s_js_is_j$.
Then we know that $R(m_1,m_2,m_3)=(n_1,n_2,n_3)$ and that
$R(m'_1,m'_2,m'_3)=(n'_1,n'_2,n'_3)$, where
$R=R_{s_1s_2s_1}^{s_2s_1s_2}$ is Lusztig's piecewise reparametrization
function associated to the canonical basis in type $A_2$ (see~\cite{lusztig2}).
Let $\Delta'_1,\Delta'_2,\Delta'_3$ be the operators on $\mathcal{B}^*$
associated to the reduced decomposition $s_js_is_j$ as in 
Definition~\ref{Delta}.
It follows from Lemma~\ref{PBWstring} that
$b_0^*=\Delta'_3\Delta'_2\Delta'_1(b^*)$, and that
$b'^*_0=\Delta'_3\Delta'_2\Delta'_1(b'^*)$. Hence, a similar argument to that
given above shows that
\begin{equation}\label{difference2}
n_0-n=r'_1-r'_2,
\end{equation}
where
\begin{eqnarray*}
r'_1 & = & n_1n'_1+n_2n'_2+n_3n'_3+n_1n'_2+n_2n'_3-n_1n'_3- \\
& & n'_1(\gamma,\alpha_i)-n'_2(\gamma,\alpha_i+\alpha_j)-
n_3(\gamma,\alpha_j)
\end{eqnarray*}
and
\begin{eqnarray*}
r'_2 & = & n'_1n_1+n'_2n_2+n'_3n_3+n'_1n_2+n'_2n_3-n'_1n_3- \\
& & n_1(\gamma',\alpha_i)-n_2(\gamma',\alpha_i+\alpha_j)-
n_3(\gamma',\alpha_j).
\end{eqnarray*}
Suppose now that $m_1>m_3$ and $m'_1<m'_3$. Then
$(n_1,n_2,n_3)=R(m_1,m_2,m_3)=(m_2,m_3,m_1+m_2-m_3)$ and
$(n'_1,n'_2,n'_3)=R(m'_1,m'_2,m'_3)=(m'_2+m'_3-m'_1,m'_1,m'_2)$.
Equating~\ref{difference1} and~\ref{difference2} we obtain:
$(m_1-m_3)(m'_3-m'_1)=0$, which is a contradiction, as each of
$m_1-m_3$ and $m'_3-m'_1$ has been assumed to be positive. A similar
argument shows that we cannot have $m_1<m_3$ and $m'_1>m'_3$.
It follows that $\mathbf{m}$ and $\mathbf{m}'$ cannot be separated by a
PBW $\ww$-wall, and the Theorem is proved.~$\Box$
\section{Properties of PBW $\ww$-linearity domains}
In this section we show that the set of $\ww$-linearity domains forms a
fan. We start with a key proposition (which we shall also need in
section 5 when we discuss $q$-commuting properties of the dual canonical
basis).
We call the connected components of the set
of PBW $\ww$-regular points PBW $\ww$-chambers, so that PBW $\ww$-linearity
domains are the closures of PBW $\ww$-chambers
(see~\cite[\S 8]{berensteinzelevinsky1}).
\begin{lemma} \label{domainlemma}
Let $\mathbf{m},\mathbf{m}'\in\mathbb{R}_{\geq 0}^N$.
Then $\mathbf{m}$ and $\mathbf{m}'$ lie in a single PBW
$\ww$-linearity domain if and only if for all $\ww'\in R(w_0)$,
$R_{\ww}^{\ww'}(\mathbf{m})$ and $R_{\ww}^{\ww'}(\mathbf{m}')$ are
weakly on the same side of all PBW $\ww'$-walls.
\end{lemma}
\proof
By definition, $\mathbf{m}$ and $\mathbf{m}'$ lie in a single PBW $\ww$-chamber
if and only if for all $\ww'\in R(w_0)$, $R_{\ww}^{\ww'}(\mathbf{m})$ and
$R_{\ww}^{\ww'}(\mathbf{m}')$ are strictly on the same side of all PBW
$\ww'$-walls.
Since PBW $\ww$-linearity domains are the closures of PBW $\ww$-chambers,
it follows from the continuity of the functions $R_{\ww}^{\ww'}$ that
if $\mathbf{m}$ and $\mathbf{m}'$ lie in a single PBW
$\ww$-linearity domain, then for all $\ww'\in R(w_0)$,
$R_{\ww}^{\ww'}(\mathbf{m})$ and $R_{\ww}^{\ww'}(\mathbf{m}')$ are
weakly on the same side of all PBW $\ww'$-walls.
Conversely, if for all $\ww'\in R(w_0)$,
$R_{\ww}^{\ww'}(\mathbf{m})$ and $R_{\ww}^{\ww'}(\mathbf{m}')$ are weakly
on the same side of all PBW $\ww'$-walls, then for all $\ww'\in R(w_0)$,
$$R_{\ww}^{\ww'}(\mathbf{m}+\mathbf{m}')=R_{\ww}^{\ww'}(\mathbf{m})+
R_{\ww}^{\ww'}(\mathbf{m}').$$
This can be proved by induction on the number of braid relations needed
to take $\ww$ to $\ww'$ (see the proof of~\cite[8.1]{berensteinzelevinsky1}).
By Proposition~\ref{linearitydomain}, this implies that $\mathbf{m}$ and
$\mathbf{m}'$ lie in a single PBW $\ww$-linearity domain.~$\Box$
\begin{proposition} \label{samedomain}
Let $\ww\in R(w_0)$. Suppose that $\mathbf{m_i}$, $1\leq i\leq k$ and
$\sum_{i=1}^k\mathbf{m}_i$ all ie in a single PBW $\ww$-linearity domain
$X$. Suppose also that
$\sum_{i=1}^k\mathbf{m}_i$ and $\mathbf{q}$ both lie in a single
PBW $\ww$-linearity domain (not necessarily the same as $X$). Then
$\mathbf{m_i}$, $i=1,2,\ldots ,k$, $\sum_{i=1}^k\mathbf{m}_i$ and
$\mathbf{q}$ all lie in a single PBW $\ww$-linearity domain.
\end{proposition}
\proof
We use Lemma~\ref{domainlemma} throughout. By the assumptions in
the Proposition, we have: \\
(i) $R_{\ww}^{\ww'}(\mathbf{m}_i)$, $1\leq i\leq k$ and
$R_{\ww}^{\ww'}(\sum_{i=1}^k\mathbf{m}_i)$
are weakly on the same side of all PBW $\ww'$-walls, and: \\
(ii) $R_{\ww}^{\ww'}(\mathbf{q})$ and
$R_{\ww}^{\ww'}(\sum_{i=1}^k\mathbf{m}_i)$
are weakly on the same side of all PBW $\ww'$-walls.
Let $H$ be a PBW $\ww'$-wall.
Firstly, let $\ww'\in R(w_0)$ be such that
$R_{\ww}^{\ww'}(\sum_{i=1}^k\mathbf{m}_i)$
does not lie on $H$.
Then $R_{\ww}^{\ww'}(\mathbf{m}_i)$, $i=1,2,\ldots ,k$,
$R_{\ww}^{\ww'}(\sum_{i=1}^k\mathbf{m}_i)$ and $R_{\ww}^{\ww'}(\mathbf{q})$
are all weakly on the same side of $H$. Next, suppose that $\ww'\in R(w_0)$
is such that $R_{\ww}^{\ww'}(\sum_{i=1}^k \mathbf{m}_i)$ does lie
on $H$. By (i) and Proposition~\ref{linearitydomain},
$R_{\ww}^{\ww'}(\mathbf{m}_i)$ lies on $H$
for $1\leq i\leq k$. Hence $R_{\ww}^{\ww'}(\mathbf{m}_i)$, $1\leq i\leq k$,
$R_{\ww}^{\ww'}(\sum_{i=1}^k\mathbf{m}_i)$ and $R_{\ww}^{\ww'}(\mathbf{q})$
are all weakly on the same side of $H$.
It follows (using Lemma~\ref{domainlemma} again) that the
$\mathbf{m_i}$, $i=1,2,\ldots ,k$, $\sum_{i=1}^k\mathbf{m}_i$ and
$\mathbf{q}$ all lie in a single PBW $\ww$-linearity domain.~$\Box$
\par
\begin{remark}\label{qcommute} \rm
If, in the assumptions and conclusion of the proposition, the property "lie in a 
single PBW $\ww$-linearity domain" is replaced by the property
"$q$-commute" (and the tuples involved are assumed to lie in
$\mathbb{Z}_{\geq 0}^N$), the result is not clear.
\end{remark}
Recall that a strongly convex polyhedral cone is a convex polyhedral cone
$C$ for which $v\in C$ implies $-v\not\in C$.
A set of strongly convex polyhedral cones is said to form a fan if
the face of every cone in the set lies in the set and if the intersection
of any two cones in the set lies again in the set.
\begin{corollary}
Let $\ww\in R(w_0)$. Then the set of PBW $\ww$-linearity domains
(together with all their faces) forms a fan in $\mathbb{R}^N$.
\end{corollary}
\proof
We first of all note that each PBW $\ww$-linearity domain is a convex
polyhedral cone (by Proposition~\ref{linearitydomain}) and thus is in fact
a strongly convex polyhedral cone as it contains no point with negative
coordinates. The same is therefore true for all of its faces.
We now show that if $C$ and $C'$ are distinct PBW $\ww$-linearity
domains, then $C\cap C'$ is a face of $C$ and of $C'$.
Suppose that $C\cap C'$ is not a face of $C$.
Then $C\cap C'\subsetneq F$, where $F$ is a face of $C$ (since $C$
and $C'$ are convex polyhedral cones).
Let $(C')^{\circ}$ denote the interior of $C'$; we will use similar
notation for the interiors of other cones.
Fix $\mathbf{p}\in F\setminus C$. Then, we can choose $\mathbf{m}\in C\cap C'$ 
such that $\mathbf{m}+\mathbf{p}\in C\cap C'$,
since $C\cap C'$ and $F$ are convex polyhedral cones.
Now, let $\mathbf{q}$ be in $(C')^{\circ}$. We have that
$\mathbf{q}$ and $\mathbf{m}+\mathbf{p}$ lie in a single PBW
$\ww$-linearity domain, and $\mathbf{m},\mathbf{p}$ and
$\mathbf{m}+\mathbf{p}$ lie in a single PBW $\ww$-linearity domain.
So by Proposition~\ref{samedomain}, $\mathbf{q}$ and $\mathbf{p}$ both
lie in a single PBW $\ww$-linearity domain. But, as
$\mathbf{q}\in (C')^{\circ}$, $C'$ is the only PBW $\ww$-linearity domain
containing $\mathbf{q}$. Since $\mathbf{p}\not\in C'$, we have a
contradiction, and we see that $C\cap C'$ is a face of $C$
(and similarly, we can see that it is a face of $C'$).
\par
We next consider the case where $C,C'$ are faces of distinct
PBW-linearity domains, $P$ and $P'$. We show that $C\cap C'$ is a face
of $C$ and of $C'$. Suppose that $C\cap C'$ is not a face of $C$.
Then, as above, $C\cap C'\subsetneq F$, where $F$ is a face of $C$.
As above, we can then choose $\mathbf{q}\in (P')^{\circ}$,
$\mathbf{p}\in F\setminus C'$ and $\mathbf{m}\in C\cap C'$ such that
$\mathbf{m}+\mathbf{p}\in C\cap C'$, since $C\cap C'$ and $F$ are
convex polyhedral cones. Note that the first two properties imply that
$\mathbf{p}\not\in P'$. This is because $\mathbf{p}$ lies in the linear span
of $C'\cap C$, so $\mathbf{p}$ lies in the linear span of $C'$. Since
$P'$ is convex, if $\mathbf{p}$ was in $P'$ and in
the linear span of its face $C'$, it would have to lie in $C'$, a
contradiction to the assumption that $\mathbf{p}\in F\setminus C'$.
So $\mathbf{q}$ and $\mathbf{m}+\mathbf{p}$ lie in a single
PBW $\ww$-linearity domain $P'$, while $\mathbf{m},\mathbf{p}$ and
$\mathbf{m}+\mathbf{p}$ all lie in a single PBW $\ww$-linearity domain
$P$. Proposition~\ref{samedomain} asserts that $\mathbf{q}$ and
$\mathbf{p}$ both lie in a single PBW $\ww$-linearity domain, but as
$\mathbf{q}\in (P')^{\circ}$, $P'$ is the only PBW $\ww$-linearity
domain containing $\mathbf{q}$, while $\mathbf{p}\not\in P'$, so we
have a contradiction. Hence, $C\cap C'$ is a face of $C$, and similarly,
we can see that it is a face of $C'$.
Finally, suppose that $C,C'$ are distinct faces of a single PBW
$\ww$-linearity domain $P$. Then $C\cap C'$ is a face of $C$, and of
$C'$, since $P$ is a polyhedral cone. We have therefore shown that
the set of PBW $\ww$-linearity domains, together with all their faces,
forms a fan in $\mathbb{R}^N$ as required.~$\Box$
\section{$q$-commuting products of quantum flag minors}
For any reduced decomposition $\tilde w=s_{i_1}\ldots s_{i_k}$ of an element $w$ 
in $W$, we define as in~\cite{caldero2}, see also~\cite{berensteinzelevinsky2}, 
the quantum flag minors 
$\Delta_{\tilde w}^*$. Roughly speaking, $\Delta_{\tilde w}^*$ is the element of 
the dual canonical basis corresponding to the extremal vector of weight 
$w\varpi_{i_k}$ in the Weyl module with highest weight $\varpi_{i_k}$. 
By~\cite{caldero1}, we have:
\begin{proposition}
Fix a reduced decomposition $\ww=(i_1,i_2,\ldots ,i_N)$ in $R(w_0)$. Set 
$\Delta_k^*=\Delta_{\tilde w_k}^*$, where $\tilde w_k=s_{i_1}\ldots s_{i_k}$, 
$1\leq k\leq N$. Then, the algebra $A_{\ww}$ generated by
the $\Delta_k^*$ is a $q$-polynomial algebra, spanned (as a space) by a part of 
the dual canonical basis.
\end{proposition}
Let $d_{\ww}$ be the form on $\mathbb{Z}^N$ defined by 
$$d_{\ww}(\mathbf{m},\mathbf{n})=\sum_{j<i}(\beta_i,\beta_j)m_in_j+\sum_im_in_i,
$$
where the $\beta_k$'s are the roots associated to $\ww$. We denote it by $d$ if 
no confusion occurs. The form $d_{\ww}$ encodes the $q$-commutations in the 
graded algebra Gr$_{\ww}(U^+)$ (see~\cite{caldero1}) associated to $\ww$. To be 
more precise, set
$$\mathbf{n}_k=\sum_{t,\,t\leq k,\,i_t=i_k}\mathbf{e}_t.$$
Then:
\begin{proposition}\label{graded}
Fix $\ww$ in $R(w_0)$. For $k$, $1\leq k\leq N$, we have 
$\Delta_k^*=B_{\ww}(\mathbf{n}_k)^*$. Moreover, for all $\mathbf{m}$, 
$\mathbf{n}$ in $\mathbb{Z}_{\geq 0}^N$:\par\noindent
\item{(i)} $q^{d(\mathbf{m},\mathbf{n})}B(\mathbf{m})^*B(\mathbf{n})^*\in
B(\mathbf{m}+\mathbf{n})^*+\sum_{\mathbf{l}\prec\mathbf{m}+\mathbf{n}}\mathbb{Z}
[q,q^{-1}]B(\mathbf{l})^*$\par\noindent
\item{(ii)} $q^{d(\mathbf{n}_k,\mathbf{m})}\Delta_k^*E(\mathbf{m})^*\in
E(\mathbf{n}_k+\mathbf{m})^*+q\mathcal{L}^*$.
\end{proposition}
\proof
The first assertion is~\cite[2.1]{caldero1}. Part (i) is a straightforward
consequence of~\cite[Corollary 2.1 (iii) and Theorem 1.2]{caldero2}. Part (ii)
is given by~\cite[Theorem 3.2]{caldero2}.~$\Box$
\par
As in~\cite{lusztig2}, we associate a reduced decomposition $\ww(Q)$
(up to commutation) of $w_0$ to each quiver $Q$ whose underlying graph is
A$_n$. Let $R(Q)\subset R(w_0)$ be the set of all reduced decompositions
associated to $Q$.
\begin{lemma}\label{flagminors} 
All quantum flag minors can be realized in the form $B_{\ww}(\mathbf{n}_k)$,
with $\ww\in R(Q)$, for some quiver $Q$ of type $A_n$.
\end{lemma}
\proof
See~\cite[4.3]{caldero2}.~$\Box$
\par
Remark that if $\ww$ is associated to a quiver orientation, then $\prec_{\ww}$ 
is the so-called degeneration ordering,~\cite{bongartz1}, and $d_{\ww}$ is 
Reineke's homological form, see~\cite{reineke1}.
As a consequence, we have:
\begin{proposition}\label{increasing}
Let $Q$ be a quiver of type A$_n$.
Fix a reduced decomposition $\ww$ in $R(Q)$. Then, for all $k$,
$1\leq k\leq N$,\par\noindent
\item{(i)} $d_{\ww}(\mathbf{n}_k,?)$ is increasing for $\prec_{\ww}$, i.e.
$\mathbf{m}\prec_{\ww}\mathbf{n}\rightarrow d_{\ww}(\mathbf{n}_k,\mathbf{m})\leq 
d_{\ww}(\mathbf{n}_k,\mathbf{n})$,\par\noindent
\item{(ii)} $q^{d(\mathbf{n}_k,\mathbf{m})}B(\mathbf{n}_k)^*B(\mathbf{m})^*\in
B(\mathbf{n}_k+\mathbf{m})^*+q\sum_{\mathbf{l}\prec\mathbf{n}_k+\mathbf{m}}\mathbb{Z}[q]B(\mathbf{l})^*$.
\end{proposition}
\proof
(i) is~\cite[Proposition 4.2]{caldero2}. (ii) follows from (i) together with 
Proposition~\ref{graded} and Proposition~\ref{ordering} (iii).~$\Box$
\par
We can now prove our main theorem.
\begin{theorem}\label{main}
Let $c$ be a $q$-commuting product of quantum flag minors and let $b$ be an 
element of $\mathcal{B}^*$ which $q$-commutes with $c$. Then, $cb$ is an element 
of $\mathcal{B}^*$ up to a power of $q$.
\end{theorem}
\proof
In the sequel, the equalities are given "up to a power of $q$".
Set $c=b_mb_{m-1}\cdots b_1$, where the $b_k$, $1\leq k\leq m$, are
$q$-commuting quantum flag minors. Using an induction on $m$, we know that $c$ 
is an element of $\mathcal{B}^*$ up to a power of $q$.\par
Now, by Lemma~\ref{flagminors}, for each $k$ we can fix a reduced
decomposition $\ww_k$ in $R(Q)$ (for some quiver $Q$ of type A$_n$)
associated to the quantum flag minor $b_k$. Let $\mathbf{p}_k$, resp.
$\mathbf{p}$, be such that $b_k=B_{\ww_1}(\mathbf{p}_k)$, resp.
$b=B_{\ww_1}(\mathbf{p})$. We have
$c\in q^{\mathbb Z}B_{\ww_1}(\sum_k\mathbf{p}_k)$. Hence, by 
Theorem~\ref{analogue}, $\mathbf{p}$ and $\sum_k\mathbf{p}_k$ are in a single 
$\ww_1$-linearity domain, since
$b$ and $c$ $q$-commute. Moreover, this theorem implies also that 
$\mathbf{p}_r$, $1\leq r\leq m$, and $\sum_{k=1}^s\mathbf{p}_k$, $1\leq s\leq m$ 
are in a single $\ww_1$-linearity domain.
By Proposition~\ref{samedomain}, $\mathbf{p}$, $\mathbf{p}_r$, 
$\sum_{k=1}^s\mathbf{p}_k$, $1\leq r,s\leq m$, are in a single $\ww_1$-linearity 
domain. Hence :
$$(*) \hskip 5 mm R_{\ww_1}^{\ww_s}(\mathbf{p}_k+\ldots + 
\mathbf{p}_1+\mathbf{p})=R_{\ww_1}^{\ww_s}(\mathbf{p}_k)+\ldots 
+R_{\ww_1}^{\ww_s}(\mathbf{p}_1)+R_{\ww_1}^{\ww_s}(\mathbf{p}),\;\;1\leq k,s\leq 
m.$$
By Proposition~\ref{increasing} (ii),
$$B_{\ww_1}(\mathbf{p}_1)^*B_{\ww_1}(\mathbf{p})^*\in
B_{\ww_1}(\mathbf{p}_1+\mathbf{p})^*+q\sum_{\mathbf{l_1}\prec_{\ww_1}\mathbf{p}_
1+\mathbf{p}}\mathbb{Z}[q]B_{\ww_1}(\mathbf{l_1})^*.$$ 
Suppose by induction on $k$ that $$b_k\ldots b_1 b\in 
B_{\ww_k}(R_{\ww_1}^{\ww_k}(\mathbf{p}_k+\ldots + 
\mathbf{p}_1+\mathbf{p}))^*+q\sum_{\mathbf{l_k}\prec_{\ww_k}R_{\ww_1}^{\ww_k}(\mathbf{p}_k+\ldots+\mathbf{p}_1+\mathbf{p})}
\mathbb{Z}[q]B_{\ww_k}(\mathbf{l_k})^*.$$
Then,
$$b_{k+1}\ldots b_1 b\in B_{\ww_{k+1}}(R_{\ww_1}^{\ww_{k+1}}(\mathbf{p}_{k+1}))^*(B_{\ww_{k}}(R_{\ww_1}^{
\ww_k}(\mathbf{p}_k+\ldots \mathbf{p}_1+\mathbf{p}))^*\hskip 15 mm$$
$$\hskip 25 
mm+q\sum_{\mathbf{l_k}\prec_{\ww_k}R_{\ww_1}^{\ww_k}(\mathbf{p}_k+\ldots+\mathbf
{p}_1+\mathbf{p})}
\mathbb{Z}[q]B_{\ww_k}(\mathbf{l_k})^*)$$
$$=B_{\ww_{k+1}}(R_{\ww_1}^{\ww_{k+1}}(\mathbf{p}_{k+1}))^*(B_{\ww_{k+1}}(R_{\ww
_1}^{\ww_{k+1}}(\mathbf{p}_k+\ldots \mathbf{p}_1+\mathbf{p}))^*\hskip 15 mm$$
$$\hskip 25 mm 
+q\sum_{\mathbf{l_k}\prec_{\ww_k}R_{\ww_1}^{\ww_k}(\mathbf{p}_k+\ldots+\mathbf{p
}_1+\mathbf{p})}
\mathbb{Z}[q]B_{\ww_{k+1}}(R_{\ww_k}^{\ww_{k+1}}(\mathbf{l_k}))^*).$$
In the previous sum, we have, by Proposition~\ref{graded} (i) and (*):
$$R_{\ww_k}^{\ww_{k+1}}(\mathbf{l_k})\prec_{\ww_{k+1}}R_{\ww_1}^{\ww_{k+1}}
(\mathbf{p}_k)+\ldots 
R_{\ww_1}^{\ww_{k+1}}(\mathbf{p}_1)+R_{\ww_1}^{\ww_{k+1}}(\mathbf{p})=
R_{\ww_1}^{\ww_{k+1}}(\mathbf{p}_k+\ldots \mathbf{p}_1+\mathbf{p}).$$
Now, by Proposition~\ref{increasing}, and the linearity property (*), we obtain 
the induction. \par
We can now show that $cb$ satisfies the properties (i) and (ii) of 
Corollary~\ref{dualclassq}. Taking $k=m$, and using Proposition~\ref{ordering} 
and the
transitivity of $\prec_{\ww_k}$, the previous results imply that,
up to a power of $q$, $cb=
E_{\ww_{m}}(R_{\ww_1}^{\ww_m}(\mathbf{p}_m+\cdots +\mathbf{p}_1+\mathbf{p}))^*+
qx$, where $x$ is a $\mathbb{Z}[q]$-linear combination of dual PBW-basis
elements $E_{\ww_m}(\mathbf{n})^*$ with
$\mathbf{n}<R_{\ww_1}^{\ww_m}(\mathbf{p}_m+\cdots +\mathbf{p}_1+\mathbf{p})$,
and therefore $cb$ satisfies property Corollary~\ref{dualclassq}(ii).
Now, Part (i) of Corollary~\ref{dualclassq} is clear because $c$ and $b$ are
$q$-commuting elements of the dual canonical basis.~$\Box$
\par
This implies that the answer to Question~\ref{question} is yes when the real
element is in the adapted algebra associated to any reduced decomposition.
\begin{corollary} Let $\ww$ in $R(w_0)$ and suppose that $b$ in $\mathcal{B}^*$ 
and $c$ in $\mathcal{B}^*\cap A_{\ww}$ $q$-commute. Then, $b$ and $c$ are 
multiplicative.
\end{corollary}

\end{document}